\documentclass[journal]{IEEEtran}
\usepackage{cite}
\usepackage{float}
\usepackage{graphicx} 
\usepackage{tikz}
\usepackage{pgfplots}
\usepackage{color}
\usepackage{xcolor}
\usepackage{amstext,epsfig,amssymb,amsbsy,amsmath}
\usepackage{mathtools}
\usepackage[output-decimal-marker={.}]{siunitx}
\usepackage{array}
\usepackage{tabularx}
\usepackage{multirow}
\usepackage{multicol}
\usepackage{setspace}
\usepackage[caption=false,font=footnotesize]{subfig}
\usepackage{wrapfig}
\usepackage{url}
\usepackage{hyperref}
\usepackage{cleveref}

\graphicspath{{Figure/}}

\expandafter\def\expandafter\UrlBreaks\expandafter{\UrlBreaks
  \do\a\do\b\do\c\do\d\do\e\do\f\do\g\do\h\do\i\do\j%
  \do\k\do\l\do\m\do\n\do\o\do\p\do\q\do\r\do\s\do\t%
  \do\u\do\v\do\w\do\x\do\y\do\z\do\A\do\B\do\C\do\D%
  \do\E\do\F\do\G\do\H\do\I\do\J\do\K\do\L\do\M\do\N%
  \do\O\do\P\do\Q\do\R\do\S\do\T\do\U\do\V\do\W\do\X%
  \do\Y\do\Z}

\newcommand\Tstrut{\rule{0pt}{2.6ex}}         
\newcommand\Bstrut{\rule[-0.9ex]{0pt}{0pt}}   

\definecolor{gray1}{rgb}{0.23529,0.23529,0.23529}%
\definecolor{gray2}{rgb}{0.49412,0.49412,0.49412}%
\definecolor{gray3}{rgb}{0.86275,0.86275,0.86275}%
\definecolor{orange1}{rgb}{0.87059,0.49020,0.00000}%
\definecolor{blue1}{rgb}{0.20392,0.30196,0.49412}%
\definecolor{blue2}{rgb}{0.72941,0.83137,0.95686}%
\definecolor{blue3}{rgb}{0.15294,0.22745,0.37255}%
\definecolor{orange2}{rgb}{1.00000,0.60000,0.00000}%
\definecolor{orange3}{rgb}{1.00000,0.89020,0.66667}%

\pgfplotsset{
    layers/my layer set/.define layer set={
        background,
        main,
        foreground
    }{},
    set layers=my layer set,
}

\begin{document}

\bstctlcite{IEEEexample:BSTcontrol}

\title{A modified version of the IEEE 3-area RTS '96 Test Case for time series analysis}%

\author{Andrea~Tosatto,
        Tilman~Weckesser,
        and~Spyros~Chatzivasileiadis,%
\thanks{A. Tosatto and S. Chatzivasileiadis are with the Technical University of Denmark, Department of Electrical Engineering, Kgs. Lyngby, Denmark (emails: \{antosat,spchatz\}@elektro.dtu.dk). T. Weckesser is with Dansk Energi, Frederiksberg C, Denmark (email: twe@danskenergi.dk).}}%

\maketitle

\vspace{-1em}
\begin{abstract}
This report describes a modified version of the IEEE 3-area RTS '96 Test Case for time series analysis. This test case was originally developed to investigate the impact of the introduction of losses in the market clearing, thus the main application of this system is DC Optimal Power Flow (OPF) studies.
A fourth area is included in the system. In each snapshot, wind power production and load consumption are modified. This allows for different import-export situations among zones, varying the prices in each zone. Moreover, several High-Voltage Direct-Current (HVDC) lines are included in the system.
\end{abstract}

\begin{IEEEkeywords}
HVDC Transmission, Time Series Analysis, Zonal Pricing Markets.
\end{IEEEkeywords}


\section{Introduction}\label{sec:1}
\IEEEPARstart{I}{n the} last decades, over 25,000 km of High-Voltage Direct-Current (HVDC) lines have been gradually integrated to the existing pan-European HVAC system. Thanks to their properties, HVDC lines facilitate the transfer of bulk power over long distances, allow the connection of asynchronous areas, represent a cost-effective solution for long-distance submarine cables and provide a useful tool for power systems stability. HVAC-HVDC interaction will be a key feature in the coming years, both for system operation and system reliability.

In this regard, there are few test systems which include meshed HVDC links. This report presents a modified version of the IEEE RTS '96 Test System \cite{1_1}, with the inclusion of several HVDC interconnectors. The intention is to develop a test system with universal characteristics, to be used as a reference for testing the impact of different evaluation techniques on diverse applications and technologies. This test system was used for the first time in \cite{1_2}, for analyzing the impact of the introduction of HVDC losses in the market clearing algorithm of zonal pricing markets. 

The test case is modified as follows:
\begin{itemize}
\item A fourth area is included.
\item Three new VSC-based HVDC links are included to connect Area 1 to the other areas.
\item Different wind farms are added in each area; different wind profiles are considered.
\item All the loads are considered elastic; three different load profiles are considered.
\item To create high and low price areas, generator costs and load utilities in the four areas are multiplied by scaling factors.
\item To create price differences over the year, load consumption and wind generation vary in each snapshot.

\end{itemize}

\section{System Topology}\label{sec:2}
\begin{figure}[!b]
\centering
\includegraphics[trim = 0.5cm 0.5cm 0.5cm 0.5cm,clip,width=0.45\textwidth]{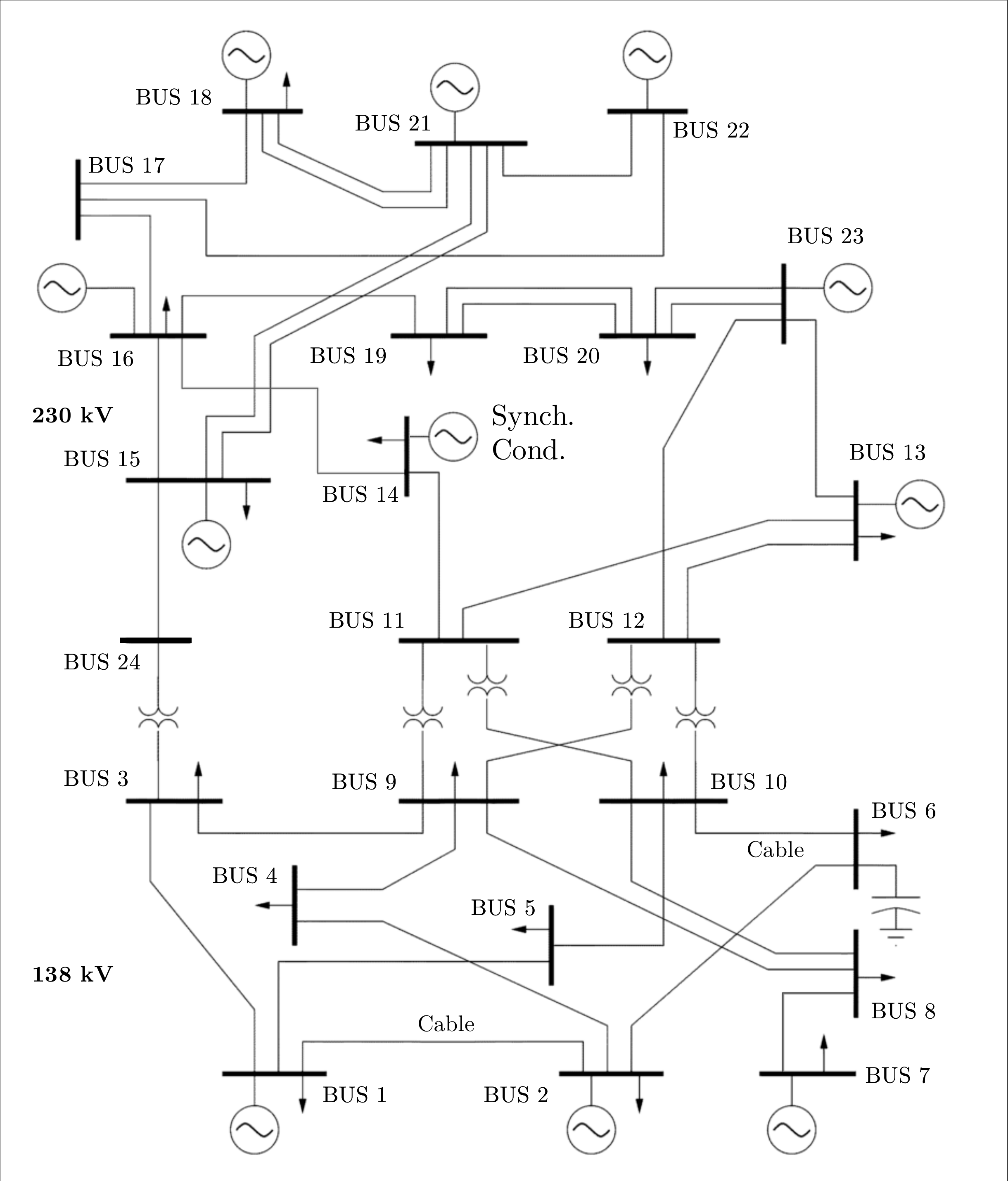}%
\vspace*{0.1em}{%
\caption{IEEE RTS 24-bus Test System.}}
\label{fig:2_onezone}
\end{figure}

\begin{figure}[h!]
\centering
\includegraphics[trim = 2.5cm 10cm 3cm 9.5cm,clip,width=1\textwidth]{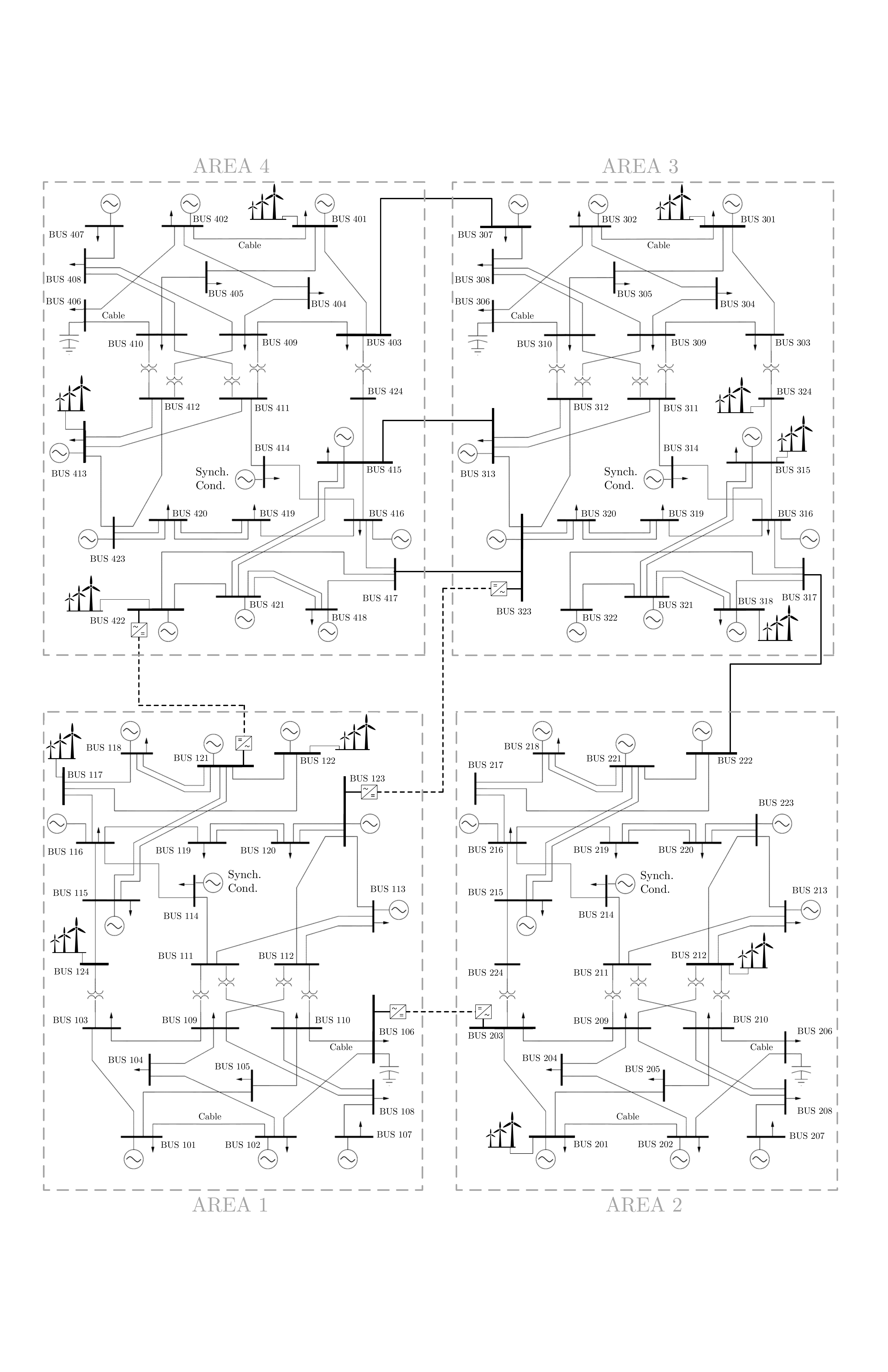}
\caption{4-area 96-bus system.}
\label{fig:2_wholesystem}
\end{figure}

This test case is intended to represent a power system where zonal pricing is applied. For this reason, the system is a multi-area power system developed by merging four different areas, each of which is the IEEE RTS 24-bus Test System shown in \figurename~\ref{fig:2_onezone}. 
The whole system is represented in \figurename~\ref{fig:2_wholesystem}. The four areas are interconnected by the following new interconnectors:
\begin{itemize}
    \item[-] 400 kV VSC-HVDC line connecting bus \# 106 and bus \# 203;
    \item[-] 250 kV VSC-HVDC line connecting bus \# 123 and bus \# 323;
    \item[-] 350 kV VSC-HVDC line connecting bus \# 121 and bus \# 422;
    \item[-] 230 kV AC line connecting bus \# 222 and bus \# 317;
    \item[-] 138 kV AC line connecting bus \# 307 and bus \# 403;
    \item[-] 230 kV AC line connecting bus \# 313 and bus \# 415;
    \item[-] 230 kV AC line connecting bus \# 323 and bus \# 417;
\end{itemize}

\clearpage
\section{Bus Data}\label{sec:3}
The bus data is obtained from the IEEE RTS '96 Test Stystem \cite{1_1}. The four areas are labeled with numbers raging from 1 to 4. The buses in each area are numbered with a preassigned numbering system (e.g., in Area 1, the 24 buses are labeled with number raging from 101 to 124, in Area 2 from 201 to 224, and so on). The slack bus is Bus \# 113 in Area 1. The bus data for Area 1 and 2 is listed in Table~\ref{tab:3_bus}. Area 3 and 4 are equal to Area 2. The bus types are the following: \begin{itemize}
    \item 1 - load bus;
    \item 2 - generator bus;
    \item 3 - slack bus.
\end{itemize}
\begin{table}[!b]
    \centering
    \caption{bus data (area 1 and 2)}
    \label{tab:3_bus}
    \footnotesize
    \begin{tabularx}{0.45\textwidth}{|*8{>{\centering\arraybackslash}X}|}
        \hline
        AREA & BUS & BUS & & & Base & V & V \Tstrut \\
        \# & \# & TYPE & GL & BL & kV & max & min \Bstrut \\
        \hline
        1	&	101	&	2	&	0	&	0	&	138	&	1.05	&	0.95	\Tstrut \\
        1	&	102	&	2	&	0	&	0	&	138	&	1.05	&	$\text{0.95 }$	\\
        1	&	103	&	1	&	0	&	0	&	138	&	1.05	&	$\text{0.95 }$	\\
        1	&	104	&	1	&	0	&	0	&	138	&	1.05	&	$\text{0.95 }$	\\
        1	&	105	&	1	&	0	&	0	&	138	&	1.05	&	$\text{0.95 }$	\\
        1	&	106	&	1	&	0	&	100	&	138	&	1.05	&	$\text{0.95 }$	\\
        1	&	107	&	2	&	0	&	0	&	138	&	1.05	&	$\text{0.95 }$	\\
        1	&	108	&	1	&	0	&	0	&	138	&	1.05	&	$\text{0.95 }$	\\
        1	&	109	&	1	&	0	&	0	&	138	&	1.05	&	$\text{0.95 }$	\\
        1	&	110	&	1	&	0	&	0	&	138	&	1.05	&	$\text{0.95 }$	\\
        1	&	111	&	1	&	0	&	0	&	230	&	1.05	&	$\text{0.95 }$	\\
        1	&	112	&	1	&	0	&	0	&	230	&	1.05	&	$\text{0.95 }$	\\
        1	&	113	&	3	&	0	&	0	&	230	&	1.05	&	$\text{0.95 }$	\\
        1	&	114	&	2	&	0	&	0	&	230	&	1.05	&	$\text{0.95 }$	\\
        1	&	115	&	2	&	0	&	0	&	230	&	1.05	&	$\text{0.95 }$	\\
        1	&	116	&	2	&	0	&	0	&	230	&	1.05	&	$\text{0.95 }$	\\
        1	&	117	&	1	&	0	&	0	&	230	&	1.05	&	$\text{0.95 }$	\\
        1	&	118	&	2	&	0	&	0	&	230	&	1.05	&	$\text{0.95 }$	\\
        1	&	119	&	1	&	0	&	0	&	230	&	1.05	&	$\text{0.95 }$	\\
        1	&	120	&	1	&	0	&	0	&	230	&	1.05	&	$\text{0.95 }$	\\
        1	&	121	&	2	&	0	&	0	&	230	&	1.05	&	$\text{0.95 }$	\\
        1	&	122	&	2	&	0	&	0	&	230	&	1.05	&	$\text{0.95 }$	\\
        1	&	123	&	2	&	0	&	0	&	230	&	1.05	&	$\text{0.95 }$	\\
        1	&	124	&	1	&	0	&	0	&	230	&	1.05	&	$\text{0.95 }$	\\
        2	&	201	&	2	&	0	&	0	&	138	&	1.05	&	$\text{0.95 }$	\\
        2	&	202	&	2	&	0	&	0	&	138	&	1.05	&	$\text{0.95 }$	\\
        2	&	203	&	1	&	0	&	0	&	138	&	1.05	&	$\text{0.95 }$	\\
        2	&	204	&	1	&	0	&	0	&	138	&	1.05	&	$\text{0.95 }$	\\
        2	&	205	&	1	&	0	&	0	&	138	&	1.05	&	$\text{0.95 }$	\\
        2	&	206	&	1	&	0	&	100	&	138	&	1.05	&	$\text{0.95 }$	\\
        2	&	207	&	2	&	0	&	0	&	138	&	1.05	&	$\text{0.95 }$	\\
        2	&	208	&	1	&	0	&	0	&	138	&	1.05	&	$\text{0.95 }$	\\
        2	&	209	&	1	&	0	&	0	&	138	&	1.05	&	$\text{0.95 }$	\\
        2	&	210	&	1	&	0	&	0	&	138	&	1.05	&	$\text{0.95 }$	\\
        2	&	211	&	1	&	0	&	0	&	230	&	1.05	&	$\text{0.95 }$	\\
        2	&	212	&	1	&	0	&	0	&	230	&	1.05	&	$\text{0.95 }$	\\
        2	&	213	&	2	&	0	&	0	&	230	&	1.05	&	$\text{0.95 }$	\\
        2	&	214	&	2	&	0	&	0	&	230	&	1.05	&	$\text{0.95 }$	\\
        2	&	215	&	2	&	0	&	0	&	230	&	1.05	&	$\text{0.95 }$	\\
        2	&	216	&	2	&	0	&	0	&	230	&	1.05	&	$\text{0.95 }$	\\
        2	&	217	&	1	&	0	&	0	&	230	&	1.05	&	$\text{0.95 }$	\\
        2	&	218	&	2	&	0	&	0	&	230	&	1.05	&	$\text{0.95 }$	\\
        2	&	219	&	1	&	0	&	0	&	230	&	1.05	&	$\text{0.95 }$	\\
        2	&	220	&	1	&	0	&	0	&	230	&	1.05	&	$\text{0.95 }$	\\
        2	&	221	&	2	&	0	&	0	&	230	&	1.05	&	$\text{0.95 }$	\\
        2	&	222	&	2	&	0	&	0	&	230	&	1.05	&	$\text{0.95 }$	\\
        2	&	223	&	2	&	0	&	0	&	230	&	1.05	&	$\text{0.95 }$	\\
        2	&	224	&	1	&	0	&	0	&	230	&	1.05	&	0.95	\Bstrut \\
        \hline
    \end{tabularx}
    \normalsize
\end{table}
GL and BL refer to the real and imaginary components of the shunt admittance to ground, Base kV to the voltage level, Vmax and Vmin to the maximum allowed voltage deviation.

\section{Elastic loads}\label{sec:4}
In each area, 17 elastic load are included. Their location and their peak values are taken from the IEEE RTS '96 Test System \cite{1_1} (peak values increased by 10\%), their utilities are linear and derived from \cite{4_1}. The load data is listed in Table~\ref{tab:4_load}. To create high and low price areas, load utilities in the four areas are multiplied by a scaling factor, respectively 1.8, 0.95, 1 and 1.1.

\begin{figure}[!h]
    \begin{tikzpicture}
        \begin{axis}[%
                width=0.16\textwidth,
                height=0.10\textheight,
                at={(1.20in,0.806in)},
                scale only axis,
                xmin=0,
                xmax=23,
                xtick={ 0,  4,  8,  12, 16, 20, 24},
                xlabel style={font=\color{white!15!black}},
                xlabel={Time (h)},
                ymin=0,
                ymax=1,
                ylabel style={font=\color{white!15!black}},
                ylabel={Load value (p.u.)},
                ylabel near ticks,
                xtick pos=left,
                ytick pos=left,
                label style={font=\footnotesize},
                every tick label/.append style={font=\footnotesize},
                axis background/.style={fill=white},
                title style={font=\bfseries\footnotesize,yshift=-1.5ex},
                title={Daily profiles}
                ]
                \addplot [color=blue1, line width=1.0pt] table[row sep=crcr]{%
                0	0.6\\
                1	0.48\\
                2	0.42\\
                3	0.38\\
                4	0.35\\
                5	0.36\\
                6	0.4\\
                7	0.46\\
                8	0.51\\
                9	0.53\\
                10	0.53\\
                11	0.54\\
                12	0.58\\
                13	0.56\\
                14	0.54\\
                15	0.55\\
                16	0.61\\
                17	0.83\\
                18	0.96\\
                19	1\\
                20	0.95\\
                21	0.9\\
                22	0.79\\
                23	0.81\\
                };
                
                \addplot [color=blue2, line width=1.0pt] table[row sep=crcr]{%
                0	0.15\\
                1	0.14\\
                2	0.11\\
                3	0.09\\
                4	0.1\\
                5	0.115\\
                6	0.22\\
                7	0.35\\
                8	0.62\\
                9	0.85\\
                10	0.9\\
                11	0.88\\
                12	0.87\\
                13	0.91\\
                14	0.95\\
                15	0.96\\
                16	0.82\\
                17	0.6\\
                18	0.42\\
                19	0.37\\
                20	0.38\\
                21	0.31\\
                22	0.19\\
                23	0.16\\
                };
                
                \addplot [color=orange2, line width=1.0pt] table[row sep=crcr]{%
                0	0.04\\
                1	0.05\\
                2	0.06\\
                3	0.07\\
                4	0.08\\
                5	0.09\\
                6	0.11\\
                7	0.14\\
                8	0.17\\
                9	0.2\\
                10	0.28\\
                11	0.38\\
                12	0.5\\
                13	0.56\\
                14	0.67\\
                15	0.85\\
                16	0.87\\
                17	0.87\\
                18	0.89\\
                19	0.85\\
                20	0.7\\
                21	0.6\\
                22	0.45\\
                23	0.22\\
                };
                
            \end{axis}
    
            \begin{axis}[%
                width=0.16\textwidth,
                height=0.10\textheight,
                at={(3in,0.806in)},
                scale only axis,
                xmin=0,
                xmax=23,
                xtick={1,5,9,13,17,21,24},
                xticklabels={{Jan},{Mar},{May},{Jul},{Sep},{Nov}},
                xticklabel style={rotate=60},
                xlabel style={font=\color{white!15!black}},
                ymin=0.4,
                ymax=1,
                ylabel style={font=\color{white!15!black}},
                ylabel={Load value (p.u.)},
                ylabel near ticks,
                xtick pos=left,
                ytick pos=left,
                label style={font=\footnotesize},
                every tick label/.append style={font=\footnotesize},
                axis background/.style={fill=white},
                title style={font=\bfseries\footnotesize,yshift=-1.5ex},
                title={Yearly profiles},
                legend columns=3,
                legend style={at={(-1.55,-0.7)}, anchor=south west, legend cell align=left, align=left, draw=black, font=\footnotesize},
                every axis legend/.append style={column sep=0.5em}
                ]
                
                \addplot [color=blue1, line width=1.0pt] table[row sep=crcr]{%
                0	0.62\\
                1	0.56\\
                2	0.523809523809524\\
                3	0.49\\
                4	0.476190476190476\\
                5	0.5\\
                6	0.571428571428571\\
                7	0.62\\
                8	0.642857142857143\\
                9	0.72\\
                10	0.857142857142857\\
                11	0.91\\
                12	0.952380952380952\\
                13	0.955\\
                14	0.976190476190476\\
                15	0.999\\
                16	1\\
                17	0.94\\
                18	0.73\\
                19	0.58\\
                20	0.55\\
                21	0.57\\
                22	0.61\\
                23	0.615\\
                };
                \addlegendentry{Residential}
                
                \addplot [color=blue2, line width=1.0pt] table[row sep=crcr]{%
                0	0.730769230769231\\
                1	0.7\\
                2	0.692307692307692\\
                3	0.71\\
                4	0.769230769230769\\
                5	0.787\\
                6	0.792307692307692\\
                7	0.82\\
                8	0.861538461538462\\
                9	0.89\\
                10	0.94\\
                11	0.96\\
                12	0.961538461538461\\
                13	0.975\\
                14	0.976923076923077\\
                15	0.992\\
                16	1\\
                17	0.93\\
                18	0.82\\
                19	0.81\\
                20	0.846153846153846\\
                21	0.82\\
                22	0.753846153846154\\
                23	0.74\\
                };
                \addlegendentry{Industrial}
                
                \addplot [color=orange2, line width=1.0pt] table[row sep=crcr]{%
                0	0.683006535947712\\
                1	0.62\\
                2	0.601307189542484\\
                3	0.63\\
                4	0.65359477124183\\
                5	0.73\\
                6	0.774183006535948\\
                7	0.79\\
                8	0.823529411764706\\
                9	0.86\\
                10	0.939542483660131\\
                11	0.95\\
                12	0.939542483660131\\
                13	0.947\\
                14	0.954575163398693\\
                15	0.98\\
                16	1\\
                17	0.9\\
                18	0.78\\
                19	0.753\\
                20	0.74\\
                21	0.68\\
                22	0.660130718954248\\
                23	0.67\\
                };
                \addlegendentry{Commercial}

            \end{axis}
        \end{tikzpicture}
    \caption{Load profiles.}
    \vspace{-0.9em}
    \label{fig:4_dprofiles}
\end{figure}
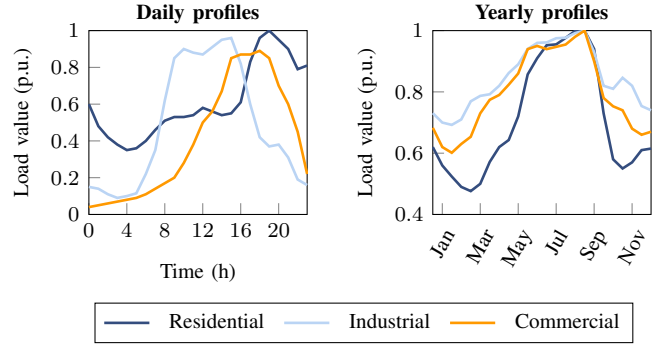

In order to create price differences over the year, the maximum consumption of loads varies according to their yearly and daily profiles (see \figurename~\ref{fig:4_dprofiles}). Three different types of load are considered: residential, industrial and commercial. The profiles are assigned to the loads using the coefficients $i_1$ and $i_2$, as shown in Table~\ref{tab:4_load}. The maximum consumption is modified as follows:

\vspace{-0.1cm}
\begin{small}
\begin{equation}\label{eq:4_load}
\begin{split}
    & D^{max} = 0.65 \cdot \textsc{res}_y \textsc{res}_d\,r\,i_1 \overline{D}^{max}\,\,+\\
    & \qquad \qquad 0.35 \cdot \textsc{com}_y \textsc{com}_d\,r\,i_1 \overline{D}^{max}\,\,+ \\
    & \qquad \qquad \textsc{ind}_y \textsc{ind}_d\,r\,i_2 \overline{D}^{max},
\end{split}
\end{equation}
\end{small}

\vspace{-0.1cm}
\noindent
where $r$ is a random numbers between 0.95 and 1.05, $i_1$ and $i_2$ are either 0 or 1 according to the load profile, $\textsc{res}_y$, $\textsc{res}_d$, $\textsc{com}_y$, $\textsc{com}_d$, $\textsc{ind}_y$ and $\textsc{ind}_d$ are the different yearly and daily coefficients, listed in Table~\ref{tab:4_profy} and \ref{tab:4_profd}, and $\overline{D}^{max}$ is the peak load, shown in Table~\ref{tab:4_load}.

\begin{table}[!b]
    \centering
    \caption{load data}
    \label{tab:4_load}
    \footnotesize
    \begin{tabularx}{0.45\textwidth}{|*2{>{\centering\arraybackslash}X} *2{c} *2{>{\centering\arraybackslash}X}|}
        \hline 
        LOAD & BUS & Peak Load & Utility  &  & \Tstrut\\
        ID & \# & \footnotesize{(MW)} & \footnotesize{(\$/MWh)} & $i_1$ & $i_2 $ \Bstrut\\
        \hline 
        D01 & 1  & 118.8 & 39.21  & 0 & 1 \Tstrut\\
        D02 & 2  & 106.7 & 35.21  & 1 & \text{$0$ } \\
        D03 & 3  & 198   & 65.35  & 1 & \text{$0$ } \\
        D04 & 4  & 81.4  & 26.89  & 0 & \text{$1$ } \\
        D05 & 5  & 78.1  & 25.79  & 1 & \text{$0$ } \\
        D06 & 6  & 149.6 & 49.42  & 1 & \text{$0$ } \\
        D07 & 7  & 137.5 & 45.42  & 0 & \text{$1$ } \\
        D08 & 8  & 188.1 & 62.11  & 1 & \text{$0$ } \\
        D09 & 9  & 192.5 & 63.56  & 0 & \text{$1$ } \\
        D10 & 10 & 214.5 & 80.82  & 1 & \text{$0$ } \\
        D11 & 13 & 291.5 & 96.97  & 0 & \text{$1$ } \\
        D12 & 14 & 213.4 & 70.45  & 1 & \text{$0$ } \\
        D13 & 15 & 348.7 & 112.07 & 0 & \text{$1$ } \\
        D14 & 16 & 110   & 36.34  & 1 & \text{$0$ } \\
        D15 & 18 & 366.3 & 135.88 & 0 & \text{$1$ } \\
        D16 & 19 & 199.1 & 65.71  & 1 & \text{$0$ } \\
        D17 & 20 & 140.8 & 46.49  & 1 & 0 \Bstrut\\
        \hline
    \end{tabularx}
    \normalsize
\end{table}

\begin{table}[!t]
    \centering
    \caption{yearly load profiles}
    \label{tab:4_profy}
    \footnotesize
    \begin{tabularx}{0.45\textwidth}{|*1{c} *3{>{\centering\arraybackslash}X}|}
        \hline 
        Period & $RES_y$ & $IND_y$ & $COM_y$ \Tstrut\Bstrut\\
        \hline 
        Jan 1 - Jan 15	&	0.620	&	0.731	&	0.683	\Tstrut \\
        Jan 16 - Jan 31	&	0.560	&	0.700	&	$\text{0.620 }$	\\
        Feb 1 - Feb 14	&	0.524	&	0.692	&	$\text{0.601 }$	\\
        Feb 15 - Feb 28	&	0.490	&	0.710	&	$\text{0.630 }$	\\
        Mar 1 - Mar 15	&	0.476	&	0.769	&	$\text{0.654 }$	\\
        Mar 16 - Mar 31	&	0.500	&	0.787	&	$\text{0.730 }$	\\
        Apr 1 - Apr 15	&	0.571	&	0.792	&	$\text{0.774 }$	\\
        Apr 16 - Apr 30	&	0.620	&	0.820	&	$\text{0.790 }$	\\
        May 1 - May 15	&	0.643	&	0.862	&	$\text{0.824 }$	\\
        May 16 - May 31	&	0.720	&	0.890	&	$\text{0.860 }$	\\
        Jun 1 - Jun 15	&	0.857	&	0.940	&	$\text{0.940 }$	\\
        Jun 16 - Jun 30	&	0.910	&	0.960	&	$\text{0.950 }$	\\
        Jul 1 - Jul 15	&	0.952	&	0.962	&	$\text{0.940 }$	\\
        Jul 16 - Jul 31	&	0.955	&	0.975	&	$\text{0.947 }$	\\
        Aug 1 - Aug 15	&	0.976	&	0.977	&	$\text{0.955 }$	\\
        Aug 16 - Aug 31	&	0.999	&	0.992	&	$\text{0.980 }$	\\
        Sep 1 - Sep 15	&	1.000	&	1.000	&	$\text{1.000 }$	\\
        Sep 16 - Sep 30	&	0.940	&	0.930	&	$\text{0.900 }$	\\
        Oct 1 - Oct 15	&	0.730	&	0.820	&	$\text{0.780 }$	\\
        Oct 16 - Oct 31	&	0.580	&	0.810	&	$\text{0.753 }$	\\
        Nov 1 - Nov 15	&	0.550	&	0.846	&	$\text{0.740 }$	\\
        Nov 16 - Nov 30	&	0.570	&	0.820	&	$\text{0.680 }$	\\
        Dec 1 - Dec 15	&	0.610	&	0.754	&	$\text{0.660 }$	\\
        Dec 16 - Dec 31	&	0.615	&	0.740	&	0.670	\Bstrut \\
        \hline
    \end{tabularx}
    \normalsize
\end{table}
\begin{table}[!h]
    \centering
    \caption{daily load profiles}
    \label{tab:4_profd}
    \footnotesize
    \begin{tabularx}{0.45\textwidth}{|*4{>{\centering\arraybackslash}X}|}
        \hline 
        Hour & $RES_d$ & $IND_d$ & $COM_d$ \Tstrut\Bstrut\\
        \hline 
        00-01	&	0.600	&	0.150	&	0.040	\Tstrut \\
        01-02	&	0.480	&	0.140	&	$\text{0.050 }$	\\
        02-03	&	0.420	&	0.110	&	$\text{0.060 }$	\\
        03-04	&	0.380	&	0.090	&	$\text{0.070 }$	\\
        04-05	&	0.350	&	0.100	&	$\text{0.080 }$	\\
        05-06	&	0.360	&	0.115	&	$\text{0.090 }$	\\
        06-07	&	0.400	&	0.220	&	$\text{0.110 }$	\\
        07-08	&	0.460	&	0.350	&	$\text{0.140 }$	\\
        08-09	&	0.510	&	0.620	&	$\text{0.170 }$	\\
        09-10	&	0.530	&	0.850	&	$\text{0.200 }$	\\
        10-11	&	0.530	&	0.900	&	$\text{0.280 }$	\\
        11-12	&	0.540	&	0.880	&	$\text{0.380 }$	\\
        12-13	&	0.580	&	0.870	&	$\text{0.500 }$	\\
        13-14	&	0.560	&	0.910	&	$\text{0.560 }$	\\
        14-15	&	0.540	&	0.950	&	$\text{0.670 }$	\\
        15-16	&	0.550	&	0.960	&	$\text{0.850 }$	\\
        16-17	&	0.610	&	0.820	&	$\text{0.870 }$	\\
        17-18	&	0.830	&	0.600	&	$\text{0.870 }$	\\
        18-19	&	0.960	&	0.420	&	$\text{0.890 }$	\\
        19-20	&	1.000	&	0.370	&	$\text{0.850 }$	\\
        20-21	&	0.950	&	0.380	&	$\text{0.700 }$	\\
        21-22	&	0.900	&	0.310	&	$\text{0.600 }$	\\
        22-23	&	0.790	&	0.190	&	$\text{0.450 }$	\\
        23-24	&	0.810	&	0.160	&	0.220	\Bstrut\\
        \hline
    \end{tabularx}
    \normalsize
\end{table}
\section{Generating units}\label{sec:5}
In each area, 33 generating units are included. Their data is taken from the IEEE RTS '96 Test System \cite{1_1} and listed in Table~\ref{tab:5_gen}. The minimum output level of all generating units has been set to zero. The cost of production of each unit is assumed to be a linear function of their output level. As for loads, these costs are multiplied by different scaling factors, respectively 0.97, 1.03, 1 and 0.99, according to the area where generators are located in. Start-up and shut-down costs, unit cylce restrictions, ramping rates and emissions have not been considered. Unit G15 is a synchronous condenser.

\begin{table}[!b]
    \vspace{-0.2cm}
    \centering
    \caption{generator data}
    \label{tab:5_gen}
    \footnotesize
    \begin{tabularx}{0.45\textwidth}{|*5{>{\centering\arraybackslash}X} {c}|}
        \hline 
        UNIT & BUS & P max  & Q min & Q max & Cost \Tstrut\\
        ID & \# & \footnotesize{(MW)} & \footnotesize{(MVAR)} & \footnotesize{(MVAR)} & \footnotesize{(\$/MWh)} \Bstrut\\
        \hline 
        G01	&	1	&	20	&	0	&	10	&	130	\Tstrut \\
        G02	&	1	&	20	&	0	&	10	&	$\text{130.00 }$ \\	
        G03	&	1	&	76	&	-25	&	30	&	$\text{16.08 }$	\\
        G04	&	1	&	76	&	-25	&	30	&	$\text{16.08 }$	\\
        G05	&	2	&	20	&	0	&	10	&	$\text{130.00 }$ \\	
        G06	&	2	&	20	&	0	&	10	&	$\text{130.00 }$ \\	
        G07	&	2	&	76	&	-25	&	30	&	$\text{16.08 }$	\\
        G08	&	2	&	76	&	-25	&	30	&	$\text{16.08 }$	\\
        G09	&	7	&	100	&	0	&	60	&	$\text{43.66 }$	\\
        G10	&	7	&	100	&	0	&	60	&	$\text{43.66 }$	\\
        G11	&	7	&	100	&	0	&	60	&	$\text{43.66 }$	\\
        G12	&	13	&	197	&	0	&	80	&	$\text{48.58 }$	\\
        G13	&	13	&	197	&	0	&	80	&	$\text{48.58 }$	\\
        G14	&	13	&	197	&	0	&	80	&	$\text{48.58 }$	\\
        G15	&	14	&	0	&	-50	&	200	&	$\text{0.00 }$	\\
        G16	&	15	&	12	&	0	&	6	&	$\text{56.56 }$	\\
        G17	&	15	&	12	&	0	&	6	&	$\text{56.56 }$	\\
        G18	&	15	&	12	&	0	&	6	&	$\text{56.56 }$	\\
        G19	&	15	&	12	&	0	&	6	&	$\text{56.56 }$	\\
        G20	&	15	&	12	&	0	&	6	&	$\text{56.56 }$	\\
        G21	&	15	&	155	&	-50	&	80	&	$\text{12.39 }$	\\
        G22	&	16	&	155	&	-50	&	80	&	$\text{12.39 }$	\\
        G23	&	18	&	400	&	-50	&	200	&	$\text{4.42 }$	\\
        G24	&	21	&	400	&	-50	&	200	&	$\text{4.42 }$	\\
        G25	&	22	&	50	&	-10	&	16	&	$\text{1E-04 }$	\\
        G26	&	22	&	50	&	-10	&	16	&	$\text{1E-04 }$	\\
        G27	&	22	&	50	&	-10	&	16	&	$\text{1E-04 }$	\\
        G28	&	22	&	50	&	-10	&	16	&	$\text{1E-04 }$	\\
        G29	&	22	&	50	&	-10	&	16	&	$\text{1E-04 }$	\\
        G30	&	22	&	50	&	-10	&	16	&	$\text{1E-04 }$	\\
        G31	&	23	&	155	&	-50	&	80	&	$\text{12.39 }$	\\
        G32	&	23	&	155	&	-50	&	80	&	$\text{12.39 }$	\\
        G33	&	23	&	350	&	-25	&	150	&	11.85 \Bstrut \\	
        \hline
    \end{tabularx}
    \normalsize
\end{table}

\begin{table}[!b]
    \centering
    \caption{wind farm data}
    \label{tab:6_wind}
    \footnotesize
    \begin{tabularx}{0.45\textwidth}{|*4{>{\centering\arraybackslash}X}|}
        \hline 
        AREA & BUS & P max & Wind \Tstrut\\
        \# & \# & \footnotesize{(MW)} & coefficient\Bstrut\\
        \hline
        1	&	117	&	113.5	&	DK1	\Tstrut\\
        1	&	122	&	56.75	&	$\text{DK1 }$	\\
        1	&	122	&	56.75	&	$\text{DK1 }$	\\
        1	&	124	&	113.5	&	DK1	\Bstrut\\
        2	&	201	&	17.03	&	DK2	\Tstrut\\
        2	&	212	&	8.51	&	$\text{DK2 }$	\\
        2	&	212	&	8.52	&	DK2	\Bstrut\\
        3	&	301	&	34.05	&	SE1	\Tstrut\\
        3	&	301	&	34.05	&	$\text{SE1 }$	\\
        3	&	315	&	68.1	&	$\text{SE1 }$	\\
        3	&	317	&	68.1	&	$\text{SE1 }$	\\
        3	&	324	&	68.1	&	SE1	\Bstrut\\
        4	&	401	&	17.25	&	SE4	\Tstrut\\
        4	&	401	&	17.25	&	$\text{SE4 }$	\\
        4	&	413	&	34.05	&	$\text{SE4 }$	\\
        4	&	422	&	34.05	&	SE4	\Bstrut\\

        \hline
    \end{tabularx}
    \normalsize
\end{table}

\section{Wind farms}\label{sec:6}
\begin{figure}[!t]
\centering
\includegraphics[trim = 2cm 16.9cm 11cm 2cm,clip,width=0.48\textwidth]{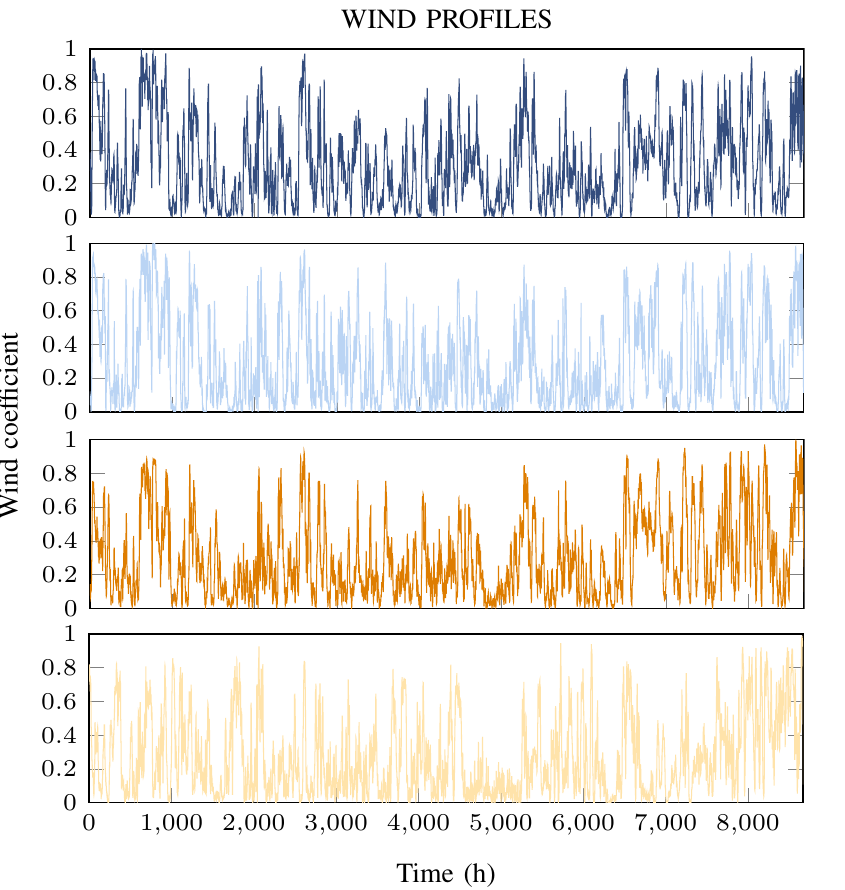}
\caption{Wind profiles. From the top: Area 1 (DK1), Area 2 (DK2), Area 3 (SE1) and Area 4 (SE4) \cite{6_1, 6_2}.}
\label{fig:6_wprof}
\end{figure}

In each area, several wind farms are included. Their data is listed in Table~\ref{tab:6_wind}. No uncertainty is considered: the output of the wind farm is known, and it varies according to the wind profile of each area. In order to have reasonable wind profiles, they are based on the wind power production of DK1, DK2 \cite{6_1}, SE1 and SE4 \cite{6_2} in 2016. The four wind profiles are depicted in \figurename~\ref{fig:6_wprof}.

\section{Transmission system}\label{sec:7}
The internal network of each area is a meshed AC grid with 31 overhead lines, 2 cables and 5 transformers (see \figurename~\ref{fig:2_onezone}). Transmission network data is obtained from the IEEE RTS '96 Test System \cite{1_1} and listed in Table~\ref{tab:7_line}. Only continuous line ratings are considered.

The four area are interconnected by the following AC and HVDC tie-lines:
\begin{itemize}
    \item[-] 400 kV VSC-HVDC line connecting bus \# 106 and bus \# 203;
    \item[-] 250 kV VSC-HVDC line connecting bus \# 123 and bus \# 323;
    \item[-] 350 kV VSC-HVDC line connecting bus \# 121 and bus \# 422;
    \item[-] 230 kV AC line connecting bus \# 222 and bus \# 317;
    \item[-] 138 kV AC line connecting bus \# 307 and bus \# 403;
    \item[-] 230 kV AC line connecting bus \# 313 and bus \# 415;
    \item[-] 230 kV AC line connecting bus \# 323 and bus \# 417;
\end{itemize}

The four AC interconnectors are adopted from \cite{1_1}, their data is listed in Table~\ref{tab:7_AC}. The base power is $\SI{100}{MW}$ and the base voltage $\SI{230}{\kilo\volt}$ (except for line 307-403, for which the base voltage is $\SI{138}{\kilo\volt}$). 

Three HVDC interconnectors are connecting Area 1 to the other areas. The base power is $\SI{100}{MW}$ and the base voltages (DC side) respectively $\SI{400}{\kilo\volt}$, $\SI{250}{\kilo\volt}$ and $\SI{350}{\kilo\volt}$. The converter station losses are represented with the generalized loss model \cite{7_1, 7_2}:
\begin{equation}\label{eq:4_1_2}
    p^{loss}_{conv} = a\,\left|I_{conv}\right|^2 + b\,\left|I_{conv}\right| + c,
\end{equation}
where $a$, $b$ and $c$ are numerical parameters reflecting the quadratic, linear and constant dependence of the losses on the line current and $\left|I^{conv}\right|$ is the magnitude of the current flowing through the converter. The data of HVDC interconnectors is derived from \cite{7_2} and listed in Table~\ref{tab:7_HVDC}.

It is worth mentioning that, since Area 1 is connected to the other areas only through HVDC lines, the whole system can represent two asynchronous areas (Area 1 and Area 2-4).

\begin{table}[!b]
    \vspace{-0.2cm}
    \centering
    \caption{line data}
    \label{tab:7_line}
    \footnotesize
    \begin{tabularx}{0.45\textwidth}{|*1{>{\centering\arraybackslash}X} *1{c} *4{>{\centering\arraybackslash}X} *2{c}|}
        \hline 
         & FROM & TO  & R & X & B & Rating & Transf \Tstrut\\
        ID & BUS & BUS & \footnotesize{(pu)} & \footnotesize{(pu)} & \footnotesize{(pu)} & \footnotesize{(MVA)} & ratio \Bstrut\\
        \hline 
        C01	&	1	&	2	&	0.003	&	0.014	&	0.461	&	175	&	0	\Tstrut\\
        L01	&	1	&	3	&	0.055	&	0.211	&	0.057	&	175	&	$\text{0 }$	\\
        L02	&	1	&	5	&	0.022	&	0.085	&	0.023	&	175	&	$\text{0 }$	\\
        L03	&	2	&	4	&	0.033	&	0.127	&	0.034	&	175	&	$\text{0 }$	\\
        L04	&	2	&	6	&	0.050	&	0.192	&	0.052	&	175	&	$\text{0 }$	\\
        L05	&	3	&	9	&	0.031	&	0.119	&	0.032	&	175	&	$\text{0 }$	\\
        T01	&	3	&	24	&	0.002	&	0.084	&	0.000	&	400	&	$\text{1.03 }$	\\
        L06	&	4	&	9	&	0.027	&	0.104	&	0.028	&	175	&	$\text{0 }$	\\
        L07	&	5	&	10	&	0.023	&	0.088	&	0.024	&	175	&	$\text{0 }$	\\
        C02	&	6	&	10	&	0.014	&	0.061	&	2.459	&	175	&	$\text{0 }$	\\
        L08	&	7	&	8	&	0.016	&	0.061	&	0.017	&	175	&	$\text{0 }$	\\
        L09	&	8	&	9	&	0.043	&	0.165	&	0.045	&	175	&	$\text{0 }$	\\
        L10	&	8	&	10	&	0.043	&	0.165	&	0.045	&	175	&	$\text{0 }$	\\
        T02	&	9	&	11	&	0.002	&	0.084	&	0.000	&	400	&	$\text{1.03 }$	\\
        T03	&	9	&	12	&	0.002	&	0.084	&	0.000	&	400	&	$\text{1.03 }$	\\
        T04	&	10	&	11	&	0.002	&	0.084	&	0.000	&	400	&	$\text{1.02 }$	\\
        T05	&	10	&	12	&	0.002	&	0.084	&	0.000	&	400	&	$\text{1.02 }$	\\
        L11	&	11	&	13	&	0.006	&	0.048	&	0.100	&	500	&	$\text{0 }$	\\
        L12	&	11	&	14	&	0.005	&	0.042	&	0.088	&	500	&	$\text{0 }$	\\
        L13	&	12	&	13	&	0.006	&	0.048	&	0.100	&	500	&	$\text{0 }$	\\
        L14	&	12	&	23	&	0.012	&	0.097	&	0.203	&	500	&	$\text{0 }$	\\
        L15	&	13	&	23	&	0.011	&	0.087	&	0.182	&	500	&	$\text{0 }$	\\
        L16	&	14	&	16	&	0.005	&	0.039	&	0.082	&	500	&	$\text{0 }$	\\
        L17	&	15	&	16	&	0.002	&	0.017	&	0.036	&	500	&	$\text{0 }$	\\
        L18	&	15	&	21	&	0.006	&	0.049	&	0.103	&	500	&	$\text{0 }$	\\
        L19	&	15	&	21	&	0.006	&	0.049	&	0.103	&	500	&	$\text{0 }$	\\
        L20	&	15	&	24	&	0.007	&	0.052	&	0.109	&	500	&	$\text{0 }$	\\
        L21	&	16	&	17	&	0.003	&	0.026	&	0.055	&	500	&	$\text{0 }$	\\
        L22	&	16	&	19	&	0.003	&	0.023	&	0.049	&	500	&	$\text{0 }$	\\
        L23	&	17	&	18	&	0.002	&	0.014	&	0.030	&	500	&	$\text{0 }$	\\
        L24	&	17	&	22	&	0.014	&	0.105	&	0.221	&	500	&	$\text{0 }$	\\
        L25	&	18	&	21	&	0.003	&	0.026	&	0.055	&	500	&	$\text{0 }$	\\
        L26	&	18	&	21	&	0.003	&	0.026	&	0.055	&	500	&	$\text{0 }$	\\
        L27	&	19	&	20	&	0.005	&	0.040	&	0.083	&	500	&	$\text{0 }$	\\
        L28	&	19	&	20	&	0.005	&	0.040	&	0.083	&	500	&	$\text{0 }$	\\
        L29	&	20	&	23	&	0.003	&	0.022	&	0.046	&	500	&	$\text{0 }$	\\
        L30	&	20	&	23	&	0.003	&	0.022	&	0.046	&	500	&	$\text{0 }$	\\
        L31	&	21	&	22	&	0.009	&	0.068	&	0.142	&	500	&	0	\Bstrut\\
        \hline
    \end{tabularx}
    \normalsize
\end{table}

\begin{table}[!b]
    \vspace{-0.2cm}
    \centering
    \caption{AC interconnector data}
    \label{tab:7_AC}
    \footnotesize
    \begin{tabularx}{0.45\textwidth}{|{c} *5{>{\centering\arraybackslash}X} {c}|}
        \hline 
        LINE & FROM & TO  & R & X & B & Rating \Tstrut\\
        ID & BUS & BUS & \footnotesize{(pu)} & \footnotesize{(pu)} & \footnotesize{(pu)} & \footnotesize{(MVA)} \Bstrut\\
        \hline 
        AC01	&	222	&	317	&	0.013	&	0.104	&	0.218	&	500	\Tstrut \\
        AC02-1	&	307	&	403	&	0.042	&	0.161	&	0.044	&	$\text{175 }$	\\
        AC02-2	&	313	&	415	&	0.010	&	0.075	&	0.158	&	$\text{500 }$	\\
        AC02-3	&	323	&	417	&	0.010	&	0.074	&	0.155	&	500	\Bstrut \\
        \hline
    \end{tabularx}
    \normalsize
\end{table}

\begin{table}[!b]
    \vspace{-0.2cm}
    \centering
    \caption{HVDC interconnector data}
    \label{tab:7_HVDC}
    \footnotesize
    \begin{tabularx}{0.45\textwidth}{|*3{c} *5{>{\centering\arraybackslash}X} {c}|}
        \hline 
        LINE & FROM & TO  & R & $a_{inv}$ & $a_{rec}$ & $b$ & $c$ & Rating \Tstrut\\
        ID & BUS & BUS & \footnotesize{(pu)} & \footnotesize{(pu)} & \footnotesize{(pu)} & \footnotesize{(pu)} & \footnotesize{(pu)} & \footnotesize{(MW)} \Bstrut\\
        \hline 
        DC01	&	106	&	203	&	.0080	&	.0005	&	.0003	&	.0007	&	.0074	&	80	\Tstrut \\
        DC02	&	123	&	323	&	.0036	&	.0056	&	.0019	&	.0013	&	.0015	&	$\text{400 }$ \\
        DC03	&	121	&	422	&	.0037	&	.0049	&	.0019	&	.0012	&	.0015	&	600	\Bstrut \\
        \hline
    \end{tabularx}
    \normalsize
\end{table}
\section{Flow-based Market Coupling}\label{sec:8}
\input{fig_FBMC.tex}

In a zonal pricing system, the network is split into price-zones in case of congestion on certain flowgates. The intra-zonal network is not included in the model, and a single price per zone is defined. The main difference between nodal and zonal pricing is that, in case of congestion, in a nodal pricing market all the nodes are subjected to different prices, while in a zonal pricing market price differences arise only among zones, with all generators and loads subjected to their zonal price \cite{4_1}. An evolution of zonal pricing is the Flow-Based Market Coupling (FBMC), which aims at coupling different independent markets. FBMC includes two clearing processes: first the energy market clearing, where a clearing price per zone is determined according to the internal power exchanges, and second, the import and export trades via the interconnections \cite{4_1}. As for zonal pricing, the intra-zonal flows are not represented in the model; in addition, cross-border lines to another zone are aggregated into a single equivalent interconnector.

Under the assumption of FBMC, all the intra-zonal nodes are aggregated into a single equivalent node. As a consequence, the 96-bus system is reduced to a 4-area system, as shown in \figurename~\ref{fig:8_FBMC}. Also, all the AC tie-lines connecting Area 3 to Area 4 are substituted by a single equivalent interconnector.  

\section{Conclusion}\label{sec:9}
This report presents a modified version of the IEEE RTS '96 Test System \cite{1_1}, with the inclusion of several HVDC interconnectors, for time series analyses. The test case is modified introducing a fourth area, three new VSC-based HVDC links, different wind farms and elastic loads. The intent is to develop a test system with universal characteristics, to be used as a reference for testing the impact of different evaluation techniques on diverse applications and technologies. The test case is designed for DC Optimal Power Flow (OPF) studies and used for the first time in \cite{1_2}.

\bibliographystyle{myIEEEtran.bst}

\end{document}